\documentclass[journal]{IEEEtran}

\usepackage{cite}

\usepackage[cmex10]{amsmath}
\usepackage{amssymb}
\usepackage{float}
\usepackage{empheq}
\usepackage{comment}
\usepackage{color}

\usepackage[super]{nth}
\usepackage[ruled,vlined]{algorithm2e}

\ifCLASSINFOpdf
   \usepackage[pdftex]{graphicx}
\else
   \usepackage[dvips]{graphicx}
\fi

\usepackage{url}
\usepackage{hyperref}

\usepackage{amsmath}
\allowdisplaybreaks

\begin{document}

\author{Andrew~W.~Rosemberg,~\IEEEmembership{}%
Alexandre~Street,~\IEEEmembership{Senior Member,~IEEE,} %
Joaquim~Dias~Garcia,~\IEEEmembership{Student Member,~IEEE,} %
Davi~M.~Vallad\~ao,~\IEEEmembership{}%
Thuener~Silva~\IEEEmembership{}%
and~Oscar~Dowson~\IEEEmembership{}%
\thanks{Andrew Rosemberg, Alexandre Street are with the Laboratory of Applied Mathematical Programming and Statistics (LAMPS PUC-Rio) and the Department of Electrical Engineering of the Pontifical Catholic University of Rio de Janeiro (PUC-Rio), Rio de Janeiro, Brazil. (email: andrewrosemberg@gmail.com, street@puc-rio.br)

Joaquim Dias Garcia is with LAMPS PUC-Rio and the Department of Electrical Engineering of PUC-Rio, Rio de Janeiro, Brazil and with PSR, Rio de Janeiro, Brazil. (email: joaquimgarcia@psr-inc.com)

Davi Vallad\~ao and Thuener Silva are with LAMPS PUC-Rio and the Department of Industrial Engineering of PUC-Rio, Rio de Janeiro Brazil. (email: davimv@puc-rio.com, thuener@gmail.com)

Oscar Dowson was with the Department of Industrial Engineering and Management Sciences, Northwestern University, USA. (email: o.dowson@gmail.com)

The work of Andrew Rosemberg was partially supported by CAPES Foundation and by FGV Energia's project P\&D ANEEL PD-09344-1703/2017. The work of Alexandre Street was partially supported by Conselho Nacional de Desenvolvimento Cient\'ifico e Tecnol\'ogico (CNPq), Funda\c{c}\~ao de Amparo \`a Pesquisa do Estado do Rio de Janeiro (FAPERJ), and by Engie Brasil Energia S.A. through project P\&D ANEEL PD-00403-0050/2020. The work of Joaquim Dias Garcia was supported in part by the Coordena\c{c}\~ao de Aperfei\c{c}oamento de Pessoal de N\'ivel Superior - Brasil (CAPES) - Finance Code 001, and by Engie Brasil Energia S.A. through project P\&D ANEEL PD-00403-0050/2020.
}
}

\title{Assessing the Cost of Network Simplifications in \\ Long-Term Hydrothermal Dispatch Planning Models}

\markboth{Accepted for publication, IEEE Transactions on Sustainable Energy, August~2021}%
{Shell \MakeLowercase{\textit{et al.}}: Bare Demo of IEEEtran.cls for IEEE Journals}


\maketitle

\begin{abstract}

The sustainable utilization of hydro energy relies on accurate estimates of the opportunity cost of the water. This value is calculated through long-term hydrothermal dispatch problems (LTHDP), and the recent literature has raised awareness about the consequences of modeling simplifications in these problems. The inaccurate representation of Kirchhoff's voltage law under the premise of a DC power flow is an example. Under a non-linear AC model, however, the LTHDP becomes intractable, and the literature lacks an accurate evaluation method of different modeling alternatives. In this paper, we extend the state-of-the-art cost-assessment framework of network approximations for LTHDP and bring relevant and practical new insights. First, we increase the quality of the assessment by using an AC power flow to simulate and compare the performance of five policies based on different network approximations. Second, we find that the tightest network relaxation (based on semidefinite programming) is not the one exhibiting the best performance. Results show that the DC {power flow} with quadratic losses approximation exhibits the lowest expected cost and inconsistency gaps. Finally, its computational burden is lower than that exhibited by the semidefinite relaxation, whereas market distortions are significantly reduced in comparison to previously published benchmarks based on DC power flow. 


\end{abstract}



\IEEEpeerreviewmaketitle


\section{Nomenclature}
\label{sec:nc}

\subsection*{Sets and Indices}

\begin{IEEEdescription}[\IEEEusemathlabelsep\IEEEsetlabelwidth{$MMM$}]
\item[$t$] Subscript to represent period $t$.


	\item[$\mathcal{N}$] Set of bus indices $n$ or $m$.
	\item[$\mathcal{N}_0$] Singleton containing reference bus index $n$ or $m$.
	\item[$\mathcal{N}_n$] Set of buses connected to bus $n$.
	\item[$\mathcal{L}$] Set of directed branches, pairs of buses $(n,m)$, where $n > m$, i.e., pairs are ordered and appear at most once.
	\item[$\mathcal{I}, \mathcal{I}_n$] Set of thermoelectric units $i$ and subset  at bus $n$.
	\item[$\mathcal{H}, \mathcal{H}_n$] Set of hydroelectric units $j$ and subset at bus $n$.
	\item[$\mathcal{H}^U_j$] Set of upstream hydroelectric units $k$ to unit $j$.
\end{IEEEdescription}

\subsection*{Constants}

\newdimen\origiwspc%
\newdimen\origiwstr%
\origiwspc=\fontdimen2\font
\origiwstr=\fontdimen3\font

\begin{IEEEdescription}[\IEEEusemathlabelsep\IEEEsetlabelwidth{$MMM$}]
    \item[$P_{it}$, $Q_{it}$] Real and reactive power bound at generator $i$.
    \item[$C_{it}$] Cost of generator $i$.
    \item[$\overline{V}_{nt}$, $\underline{V}_{nt}$] Upper/lower magnitude voltage limits at bus $n$.
    \item[$D_{nt}$,$C^{\delta}_{nt}$] Real power load and deficit cost at bus $n$.
    \item[$Y^s_{nt}$, $Y^{sq}_{nt}$] Shunt real and imaginary admittance at bus $n$.
    \item[$G^c_{(n,m)t}$] Real part of pi-section parameters at branch $(n,m)$.
    \item[$B^c_{(n,m)t}$] Imaginary part of pi-section parameters at branch $(n,m)$.
    \item[$F_{(n,m)t}$] Apparent power limit at branch $(n,m)$.
    \item[$R_{(n,m)t}$] Resistance at branch $(n,m)$.
    \item[$X_{(n,m)t}$] Reactance at branch $(n,m)$.
    \item[$\Upsilon_{jt}$] Volume limit at hydroelectric units $j$.
    \item[${\nu}_{jt-1}$] Initial volume at hydroelectric units $j$.
    \item[$A_{jt}$] Possible inflows at hydroelectric units $j$.
    \item[$U_{jt}$] Outflow limit at hydroelectric units $j$.
    \item[$\rho_{jt}$] Production factor at hydroelectric units $j$.
\end{IEEEdescription}

\subsection*{Decision Variables}

\begin{IEEEdescription}[\IEEEusemathlabelsep\IEEEsetlabelwidth{$MMM$}]
    \item[$p_{it}$, $q_{it}$] Real and reactive power dispatch at generator $i$.
    \item[$v_{nt}$, $\theta_{nt}$] Complex voltage and Phase angle at bus $n$.
    \item[$\delta_{nt}$, $\ell_{nt}$] Real power deficit and loss\footnote{Active power loss approximation in the DCLL model, and bus shunt loss in the AC, SDP, and SOCP models.} at bus $n$.
    \item[$u_{jt}$, $s_{jt}$] Outflow and Spillage at hydroelectric units $j$.
    \item[$\nu_{jt}$] Volume at hydroelectric units $j$.
    \item[$f_{(n,m)t}$] Real power flow at branch $(n,m)$.
    \item[$f^q_{(n,m)t}$] Reactive power flow at branch $(n,m)$.
\end{IEEEdescription}

\subsection*{Operators}

\begin{IEEEdescription}[\IEEEusemathlabelsep\IEEEsetlabelwidth{$MMM$}]
	\item[$(\cdot)^{*}$] Complex conjugate.
	\item[$\angle(\cdot)$, $|\cdot|$] Angle and magnitude of a complex number.
\end{IEEEdescription}



\

\

\

\section{Introduction}
\label{sec:int}

\IEEEPARstart{A}{ssessing} the value of systems' scarce resources is a key activity of power system operators, especially those in charge of hydrothermal power systems \cite{pereira1991multi}. Such assessment aims at evaluating implicit opportunity costs of relevant systems' resources, such as water, which can be stored and used in future periods to prevent the expensive thermal generation or load shedding. Due to their fast response, hydroelectric power plants also play a crucial role in paving the way for the sustainable integration and utilization of other renewable sources such as wind and solar. In addition, due to the scale of installed generation, hydropower is one of the most relevant renewable resources in the world \cite{EIA2019}. Therefore, an accurate assessment of the water value is crucial for ensuring the long-term sustainability of energy systems \cite{street2020assessing}.

Opportunity costs of systems' relevant resources, such as water in long-term hydrothermal systems, are generally estimated by solving a long-term hydrothermal dispatch problem (LTHDP). This problem is the core of operative planning studies and is formulated as a multistage stochastic linear model \cite{pereira1991multi,Street2017}. The state-of-the-art solution method for such large-scale optimization problems is the Stochastic Dual Dynamic Programming (SDDP) algorithm of \cite{pereira1991multi}. The SDDP algorithm iteratively approximates the future cost of operation, also known as the cost-to-go function, as a piecewise linear convex function of the amount of water stored in the reservoirs. The computational efficiency of SDDP notwithstanding, the algorithm relies on strong assumptions such as convexity of the dispatch problem defining the cost-to-go function. This is a consistent limitation in the literature on the subject (we refer to \cite{Helseth2016}, where an agent sells reserve capacity in a competitive market, and \cite{Street2017}, where the $n-K$ security criterion is incorporated within the centralized cost minimization framework). Some non-convexity can be modeled by including binary variables \cite{Zou2019,Martin2019}, but this introduces associated computational challenges.

Despite the recent advances in the SDDP algorithm, the convexity requirement and large problem-size of long-term hydrothermal planning studies necessitates modeling simplifications to form a tractable multistage stochastic programming problem. However, the operative policy is usually implemented by a more detailed short-term model coupled with the opportunity cost assessment of the simplified long-term operation planning model. Consequently, the implemented operative decisions may differ significantly from the planned ones even if the same scenario takes place, and the implementation model differs from the planning model solely on one feature such as the network model \cite{7859406}. Such inconsistency can be understood as a model-misspecification risk. As demonstrated in \cite{7859406}, statistically significant time-inconsistency gaps can be induced when Kirchhoff's Voltage Law (KVL) and the $n-1$ security criterion are disregarded under the DC power flow approximation. Other relevant simplifications widely adopted in the LTHDP literature are: linear hydro production functions \cite{vitor2019} and the hazard-decision approximation of the information-revelation process of inflows \cite{street2020assessing}. 

In this work, we are interested in making a controlled study to isolate the effects of the network simplifications (\emph{ceteris paribus}). Therefore, we assume the typical LTHDP with all its customary simplifications but the network one. In this way, we can measure the effects of considering each network model assuming that all other simplifications and assumptions remain fixed.

The studies in \cite{7859406} assume a DC power flow model as the reference model for the assessment of the time-inconsistency gap when planning the system operation using the transportation network flow approximation model (NFA). However, as reported in the technical literature (see \cite{low2014convex} and \cite{dvijotham2016error}), DC power flow approximation still produces significant discrepancies and the actual operation is better represented by an AC power flow model. {In this context, conic relaxations such as the semidefinite programming (SDP) (proposed in \cite{6345272} and further studied in \cite{lavaei2013}), and the second-order cone (SOC) (proposed in \cite{1664986} and further studied in \cite{kocuk2016strong}) formulations were proposed and studied to improve the state-of-the-art tractable approximations of the AC power flow. The tractability issue of AC optimal power flow (AC OPF) models are due to its nonconvexity, which prevents the application of many techniques relying on such property (we refer to \cite{Sun2021} for a recent effort to provide a convergent algorithm to solve AC OPF).} Thus, in this work, we provide novel results on the impact of network simplifications that extends previously reported works on two fronts: 1) we study the performance of five network approximation models used in the literature, and 2) we consider the more accurate AC power-flow model as the reference model used to assess the performance of each approximation in terms of total cost, inconsistency gap, and market distortions. These extensions allow us to isolate and identify the pros and cons of different network models in terms of their approximation quality and computational burden.

To achieve the aforementioned goals, a new, reliable, and flexible implementation of the SDDP algorithm and the inconsistency-assessment fast algorithm (first proposed in \cite{7859406} to assess the inconsistency due to modeling simplifications) was needed. Unfortunately, such a tool was not available in the open literature. Therefore, an open-source SDDP tool with these features, \textit{HydroPowerModels}, was developed and made available in \cite{Rosemberg2020}. This tool enables not only the assessment of the results found in this work, but also allows researchers in the electrical industry to test new ideas, leveraging state-of-the-art solution methods and mathematical formulations for the LTHDP. Therefore, the public open-source package can be seen as a relevant side contribution of this work.

The objective of this work is to present an open-source computational framework
for testing the operative and economic impact of modeling simplifications over the network power-flow in hydrothermal power systems. Among the myriad of formulations available in the package, 
we focused on assessing the cost and operative performance of the following model approximations: NFA, currently in use by the Brazilian system operator \cite{7859406}; SOC \cite{1664986}; SDP \cite{6345272}; the DC power-flow approximation (DC) \cite{Stott2009_DC}, and the DC with line-loss power-flow approximation (DCLL) \cite{dcll_coffrin}. All the previously mentioned formulations are tested as approximations for the network model in the planning stage of the LTHDP, where the cost-to-go function is built through the SDDP algorithm. Then, we evaluate each approximation by simulating the system's operation, which minimizes costs under AC power-flow constraints (see \cite{carpentier1962contribution} and \cite{Cain2012}). 

In this context, the main contributions of this work are twofold:

\begin{enumerate}
    
    \item Extend {and generalize} the time-inconsistency cost assessment technique proposed in \cite{7859406} and make it available through an open-source computational tool. { The proposed technique used to quantify the impact of network simplifications in the context of long-term dispatch planning is, for the first time, carried out using the more accurate AC power-flow as the reference/implementation model. Based on that,} we are capable of comparing, with higher precision, the performance of five network approximations (NFA, SOC, SDP, DC, and DCLL), {whose market distortions, inconsistency gap, and long-run cost efficiency have never been studied before in the context of long-term hydrothermal dispatch planning.} Within the proposed computational framework, we can isolate the benefit of each network approximation and evaluate the tradeoff between computational burden and solution quality. 
    
    \item We identify new insights about the detrimental effect of network simplifications while identifying the tradeoff between the computational burden and solution quality for each network model. More specifically, { within the limitations of our case study,} we show that the SOC approximation induces a planning policy with poor performance and dispatch distortion when the system is fully meshed (with many cicles) and with good performance when the system is approximately radial. Under mild conditions, this is always obtained with a reasonable computational burden. The planning policy based on the SDP relaxation, the tightest convex relaxation, exhibits low market distortions but comes with a high computational burden. Finally, the DCLL approximation exhibits the lowest system cost when evaluated with the full AC model, the lowest inconsistency gap, and its computational burden is lower than that presented by the convex conic relaxations (SOC and SDP). Thus, we find this modeling approximation, based on convex quadratic programming, as an interesting and promising alternative to reduce market distortions at a reasonable computational burden. To the best of our knowledge, this is the first time that long-term hydrothermal dispatch policies based on conic relaxations and DCLL network approximations are compared. 

    
\end{enumerate}

The remainder of this work is organized as follows. Section \ref{sec:tb} presents the theoretical background for LTHDP and time-inconsistency. Section \ref{sec:pf} explains the five network approximation models derived from the AC power-flow model and presents the open-source tool to support studies of LTHDP. Section \ref{sec:cases} presents case studies. Relevant conclusions are drawn in Section \ref{sec:conclusion}.

\vspace{-0.1cm}
\section{Theoretical Background}
\label{sec:tb}
The LTHDP is modeled as a large-scale multistage stochastic optimization problem and generally solved by SDDP \cite{pereira1991multi,philpott2013solving,7859406,Street2017,street2020assessing,shapiro2019,vitor2019}. The goal is to find the optimal operation planning policy of a hydrothermal power system to meet demand throughout a long-run planning horizon. For this purpose, the policy must take into account the best use of water to ensure power balance across the network during all periods under the uncertainty of inflows, demand, fuel costs, etc. Following previously reported literature on the LTHDP, in this work we consider only the inflow uncertainty. 
\vspace{-0.3cm}
\subsection{The dispatch model} \label{subsec:tb1}
The complete mathematical formulation of the LTHDP in one stage, will be represented by the function $\mathcal{Q}_t(\nu_{t-1},\omega_t) =$
\begin{subequations} \label{OPF}
    \begin{align}
        \hspace{-1cm}\underset{\mathbf{x}_t}{\mbox{min: }} & \sum_{i \in \mathcal{I}} C_{it} \, p_{it} + \sum_{n \in \mathcal{N}} C^{\delta}_{nt} \, \delta_{nt} + \mathbb{E}[\mathcal{Q}_{t+1}(\nu_{t}, \omega_{t+1})] \label{eq_objective}\\
        \mbox{s.t.: } & \nonumber \\
        & \sum_{\substack{i \in I_n}} p_{it} + \sum_{\substack{j \in \mathcal{H}_n}} u_{jt} \, \rho_{jt} -\sum_{\substack{m \in \mathcal{N}_n}} f_{(n,m)t} - \ell_{nt} + \nonumber \\
        & \qquad  =  D_{nt} - \delta_{nt}, \;\; \forall n\in \mathcal{N} \label{eq_kcl_shunt} \\
        & \nu_{jt} + u_{jt} + s_{jt} = \nu_{j,t-1} + A_{j,t}(\omega_t) + \nonumber \\ 
        & \qquad  \sum_{\substack{k \in \mathcal{H}^U_j}} u_{kt} + \sum_{\substack{k \in \mathcal{H}^{S}_j}} s_{kt}, \;\; \forall j \in \mathcal{H} \label{eq_hydro_balance} \\
        & |f_{(n,m)t}| \leq F_{(n,m)t},\; |f_{(m,n)t}| \leq F_{(n,m)t} \;\; \forall (n,m) \in \mathcal{L} \label{eq_thermal_limit}\\
        & 0 \leq p_{it} \leq P_{it} \;\; \forall i \in I  \label{eq_gen_bounds}\\
        & 0 \leq \nu_{jt} \leq \Upsilon_{jt} \;\; \forall j \in \mathcal{H} \label{eq_volume_limit}\\
        & 0 \leq u_{jt} \leq U_{jt} \;\; \forall j \in \mathcal{H} \label{eq_outflow_limit} \\
        & \ell_{nt} \geq 0 \;\; \forall n \in \mathcal{N} \label{eq_loss_opf} \\
        & \delta_{nt} \geq 0 \;\; \forall n \in \mathcal{N} \label{eq_def_opf} \\
        & \mathbf{x}_t \in \mathcal{X}_t. \label{eq_XT}
    \end{align}
\end{subequations}
$\mathbf{x}_t$ is the stacked vector of all decision variables, $\mathbf{x}_t = [p_t, f_t, u_t, s_t, \nu_t, \delta_t, \ell_t]^T$, and $\mathcal{X}_t$ represents different power-flow constraints related to the associated network formulation. It is important to mention that in $\mathcal{X}_t$, many relevant variables such as voltage levels and reactive power not appearing in model \eqref{OPF} can be considered depending on the studied network model. In this sense, \eqref{OPF} can be seen as a general model that can be adapted to consider different power flow models within our assessment framework as will be further explained.

The objective function \eqref{eq_objective} is to minimize the sum of immediate costs represented by the costs of active power generation and the cost of energy supply deficit, and future costs represented by the cost-to-go function $\mathbb{E}[\mathcal{Q}_{t+1}(\nu_{t}, \omega_{t+1})]$.
Constraints \eqref{eq_kcl_shunt} implement Kirchhoff’s Current Law (KCL), i.e., power balance at each node. Deficit variables guarantee feasibility in case of a lack of power availability.
The water mass balance equation is implemented in \eqref{eq_hydro_balance}, where the water stored in a reservoir should equal the water previously stored plus the incoming inflows and water discharged from upstream reservoirs, minus the portion used to generate energy and the one spilled away. Constraints \eqref{eq_thermal_limit} to \eqref{eq_def_opf} represent physical limits on variables.
Finally, expression \eqref{eq_XT} represents a feasibility set used to model all network models considered in this work. The characterization of this set to represent the different network models is made in section \ref{sec:pf}.

Simplified Optimal Power Flow (OPF) network formulations, $\mathcal{X}_t$, are a necessary condition to achieve a tractable model that is also compatible with efficient solution methods such as SDDP. Planning agents use simplified models to compute, generally in monthly basis, typical operating points and the associated reservoir levels for long-term horizons in order to assess cost-to-go functions. This is consistent with the state-of-the-art literature and practice (\cite{pereira1991multi,shapiro2011analysis,philpott2013solving,shapiro2013risk,Helseth2016, Street2017, street2020assessing}). For implementing operating decisions, however, independent systems operators (ISO) seek feasible dispatches complying with more detailed network models. This is done by coupling the simplified view of the systems' future operation, implicitly considered in the cost-to-go function, with more realistic network formulations. Unfortunately, in practice, the simplifications considered in the network model to evaluate cost-to-go functions are significantly optimistic compared to the representation needed to ensure feasibility. This optimistic bias leads to the implementation of expensive (sub-optimal) time-inconsistent policies \cite{7859406}.

For instance, in the Brazilian case {(see \cite{maceiral2018twenty}), the national ISO applies} the SDDP algorithm using the NFA model representation to estimate the cost-to-go function, which is then used as input in a second model, with a more accurate representation of the network, to define decisions to be implemented. Then, in the next period, the state is updated with the actual reservoir levels and the same process is repeated again. In Chile, conversely, a similar process is carried out, but the SDDP is implemented with a DCLL model (see page 75 of \cite{PLPChile}). This illustrates how the network representation in the operational planning stage varies from system to system according to the perception of the impact each representation may bring. Notwithstanding, this rolling-horizon operating scheme potentially produces time-inconsistent policies in which implemented decisions deviate from those obtained in the planning stage embedded in the cost-to-go function. In this context, hybrid and inconsistently implemented policies may produce decisions that can be far from optimal -- for both the planning problem and for the true problem based on the detailed network model. 
\vspace{-0.2cm}
\subsection{Evaluation process} \label{subsec:tb2}
To evaluate the performance of the time-inconsistent policies induced by network simplifications, we extend the idea of the fast algorithm proposed in \cite{7859406}. Roughly, this algorithm allows us to estimate simplified cost-to-go functions and simulate the cost of planning-implementation processes where two different models are used, which is the case of a hydrothermal power system. As described in subsection \ref{subsec:tb1}, two different models for the network are used to operate a hydrothermal power system: 1) $\mathcal{X}_t^{plan}$, which is applied in the planning stage, where the water values are estimated by cost-to-go functions through backward and forward iterations in an SDDP fashion {(further, in the case study section, we test as planning model the NFA, DC, DCLL, SOCP, and SDP)}; 2) $\mathcal{X}_t^{imp}$, which is used to obtain implementable dispatch decisions using the cost-to-go functions estimated with the previously described planning model (in this work, we use the AC power flow model as the implementation model to test all aforementioned approximations as planning models). 

Aiming to isolate the simplification effect of a given network model, in this work, the difference between the planning model and the implementation model is defined by which model is considered in $\mathcal{X}_t$. Hence, it is worth highlighting that the only result of the planning stage used in the implementation phase of period $t$ is the cost-to-go function, $\mathbb{E}[\mathcal{Q}_{t+1}^{plan}]$, which was build upon the assumption that $\mathcal{X}_t = \mathcal{X}_t^{plan}$. {Therefore, due to the inexactness of some of the planing models that will be studied in this work, it is possible (and likely) that some operating points obtained in the planning stage turn out to be infeasible in reality (AC power flow). In practice, however, in these cases, other operating points (some of them, relying on higher dispatch costs or load curtailment) are obtained with the more realistic AC network model in the implementation step. This is in line with the objective of this work, namely, to measure the quality of implemented decision (using the AC network model) when relying on different network simplifications in the planning stage.}

\vspace{-0.2cm}
\begin{figure}[h!]
\centering
{\includegraphics[width=0.75\linewidth]{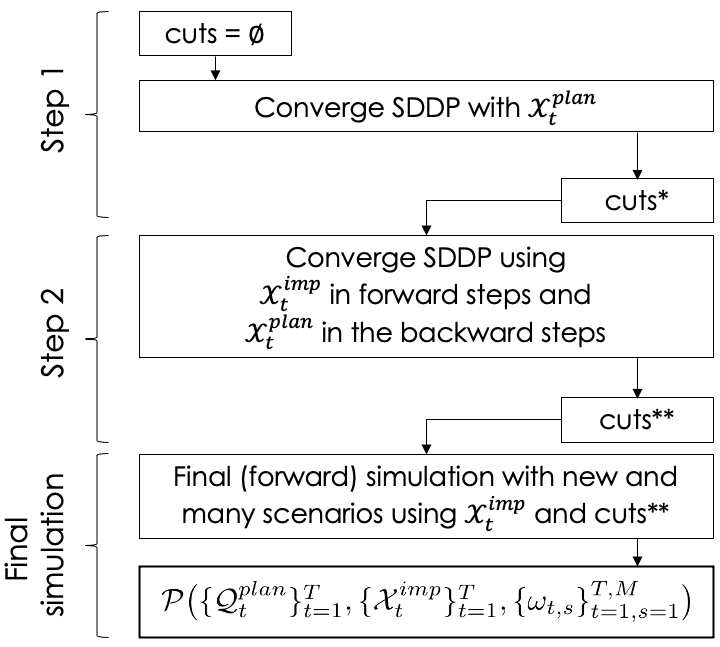}}
\vspace{-0.15cm}
\caption{{Flowchart of the updated fast algorithm (see \cite{7859406}) used in this work.}}
\label{fig:FastAlgo}
\end{figure}

To simulate the decisions that should be implemented based on this inconsistent planning-implementation process, two algorithmic steps and a final simulation (depicted in Figure \ref{fig:FastAlgo}) are carried out as follows. In the first algorithmic step {(Step 1 in Figure \ref{fig:FastAlgo})}, the SDDP with the planning model, $\mathcal{X}_t^{plan}$, is converged in the traditional way \cite{pereira1991multi,philpott2013solving,vitor2019} { and a set of cuts that approximate the cost-to-go functions of each stage are obtained (marked as cuts* in Figure \ref{fig:FastAlgo})}. Then, in the second algorithmic step {(Step 2 in Figure \ref{fig:FastAlgo})}, the cost-to-go functions obtained with the planning model, $\mathbb{E}[\mathcal{Q}_{t+1}^{plan}]$, are further approximated with new cuts {(marked as cuts** in Figure \ref{fig:FastAlgo})} through the application of an updated version of the fast algorithm described in \cite{7859406}. To do that, we initialize a modified SDDP method with the first-step previously converged cost-to-go function approximations {(cuts*)}. Then, we execute forward and backward iterations with different network models. In the forward iteration, we switch the network model to the implementation one, $\mathcal{X}_t^{imp}$, which in this work will be the AC power flow as will be further introduced. This allows us to add new cuts on more realistic points that may have been neglected by the previous step, where the cost-to-go functions were approximated with the planning model in both forward and backward iterations of the SDDP. The stopping criterion follows the upper and lower bound convergence as described in \cite{7859406} and \cite{Street2017}. After reaching the stopping criterion, a final simulation step, with many scenarios, is carried out using \eqref{OPF} with the improved cost-to-go functions, $\mathbb{E}[\mathcal{Q}_{t+1}^{plan}]$ for every $t$, and $\mathcal{X}_t^{imp}$.

Compared to \cite{7859406}, the combination of the three steps in a unified algorithm constitutes a novel and practical improvement (with reduced computational time) in the cost assessment of inconsistent policies. Furthermore, the second step is key to ensure we converge the cost-to-go function at points that are likely to be visited by the system operator when using the implementation model $\mathcal{X}^{imp}$. Note that forward passes using $\mathcal{X}^{plan}$ in Step 1 may follow trajectories that are likely to be infeasible under $\mathcal{X}^{imp}$.
Therefore, while in the first step we find the cuts generally obtained in the planning stage, the second step improves the cost-to-go function with new cuts in more realistic operating points as if the forward iteration were simulating future implementation steps. It is worth noting that, in general, good local approximations are available for the second-stage cost-to-go function.

To quantify the time-inconsistency gap, i.e., the modeling risk, we compare the cost of the actual implemented inconsistent policy, $\mathcal{P} \big(\{\mathcal{Q}_t^{plan}\}^T_{t=1}, \{\mathcal{X}_t^{imp}\}^T_{t=1}, \{\omega_{t,s}\}^{T,M}_{t=1,s=1} \big)$, with the cost of the respective reference planning policy, $\mathcal{P} \big(\{\mathcal{Q}_t^{plan}\}^T_{t=1}, \{\mathcal{X}_t^{plan}\}^T_{t=1}, \{\omega_{t,s}\}^{T,M}_{t=1,s=1} \big)$. $\mathcal{P}$ represents all the simulation results for a policy using a given set of cost-to-go functions, a given set describing the network model, and a given set of inflow scenarios.
Hence, the time-inconsistency gap measures the hidden cost of neglecting the constraints of the implementation problem in the planning phase. It can be seen as the operator's regret with respect to its planning expectations. Consequently, it allows ISOs to detect and quantify the impact of the inconsistencies induced by a given simplification without the need of simulating the full, and possibly currently intractable, policy based on the more complex network model such as the AC power flow.

\vspace{-0.2cm}
\section{Network Formulations}
\label{sec:pf}

The network constraints set $\mathcal{X}_t$ can assume many forms. We present a non-standard description of power flow constraints that is key to contrast and compare the different formulations considered in this work. 
We start with the most detailed network model, the non-linear AC formulation, which will be used as the evaluation model to evaluate the quality of all subsequent formulations. In this case, the network model is the following:\\

$\mathcal{X}_{AC,t}$ $=$ \Big\{ $p_t, f_t, u_t, s_t, \nu_t, \delta_t, \ell_t \, $\Big\vert$ \, \exists \; q_t,v_t,\mathbf{W}_t,f^q_t$ :
\begin{subequations} \label{AC_W_PF}
    \begin{align}
        & \mathrm{w}_{nmt} = v_{nt}v_{mt}  \;\; \forall(n,m) \in \mathcal{L} \label{eq_ac_w}\\
        & \underline{V}_{nt}^2 \leq \mathrm{w}_{nnt} \leq \overline{V}_{nt}^2 \;\; \forall n \in \mathcal{N} \label{eq_voltage_acw_bounds}\\
        & (f_{(n,m)t})^2 + (f^q_{(n,m)t})^2 \leq F_{(n,m)t}^2 \;\; \forall (n,m) \in \mathcal{L} \label{eq_thermal_limit_acw} \\
        & \sum_{\substack{i \in I_n}} q_{it} - Y^{sq}_{nt} \mathrm{w}_{nnt} -\sum_{\substack{m\in \mathcal{N}_n}} f^q_{(n,m)t} = 0 \;\; \forall n\in \mathcal{N} \label{eq_kcl_shunt_q_acw} \\
        & \ell_{nt} = Y^{s}_{nt} \mathrm{w}_{nnt} \;\; \forall n \in \mathcal{N} \label{eq_loss_acw}\\      
        & f_{(n,m)t} = \left( G_{(n,m)t} + G_{(n,m)t}^c \right) \mathrm{w}_{nnt} + \nonumber \\
        & \ \ -G_{(n,m)t} \mathrm{w}_{nmt}^{\Re} -B_{(n,m)t} \mathrm{w}_{nmt}^{\Im}  \;\;  \forall (n,m)\in \mathcal{L} \label{eq_power_acw}\\
        & f_{(m,n)t} = \left( G_{(m,n)t} + G_{(m,n)t}^c \right) \mathrm{w}_{mmt} +  \nonumber \\
        & \ \ -G_{(m,n)t} \mathrm{w}_{mnt}^{\Re} -B_{(m,n)t}\mathrm{w}_{mnt}^{\Im}  \;\;  \forall (n,m)\in \mathcal{L} \label{eq_power_acw_2}\\
        & f_{(n,m)t}^q = -\left( B_{(n,m)t} + B_{(n,m)t}^c\right) \mathrm{w}_{nnt} + \nonumber \\
        & \ \ +B_{(n,m)t} \mathrm{w}_{nmt}^{\Re}  -G_{(n,m)t} \mathrm{w}_{nmt}^{\Im}   \;\;  \forall (n,m)\in \mathcal{L} \label{eq_power_acw_q} \\
        & f_{(m,n)t}^q = -\left( B_{(m,n)t} + B_{(m,n)t}^c\right) \mathrm{w}_{mmt} + \nonumber\\
        & \ \ +B_{(m,n)t} \mathrm{w}_{mnt}^{\Re} - G_{(m,n)t} \mathrm{w}_{mnt}^{\Im}   \;\;  \forall (n,m)\in \mathcal{L} \label{eq_power_acw_q_2}\\
        & -Q_{it} \leq q_{it} \leq Q_{it} \;\; \forall i \in I \Big\},\label{eq_gen_q_bounds_acw}
    \end{align}
\end{subequations}
where $\mathrm{w}_{nmt}^{\Re}$ represents the real part of the $\mathrm{w}_{nmt}$ variable, and, $\mathrm{w}_{nmt}^{\Im}$, the imaginary part.

This formulation is actually an equivalent model to the more classical formulation \cite{carpentier1962contribution} developed in the attempt to provide better relaxations, \cite{BAI2008383,1664986,coffrin2015qc}. 
It uses an auxiliary variable $\mathrm{w}_{nmt}$ to represent the product of the voltage from buses $n$ and $m$, i.e., $\mathrm{w}_{nmt} = v_{nt}v_{mt}$, in constraint \eqref{eq_ac_w}. In this model, the magnitude of the complex voltage is bounded in constraint \eqref{eq_voltage_acw_bounds}, and the apparent power is limited in constraint \eqref{eq_thermal_limit_acw}. Constraints \eqref{eq_kcl_shunt_q_acw} models the reactive power flow and constraint \eqref{eq_loss_acw} represents the real part of the bus shunt loss that appears in constraint \eqref{eq_kcl_shunt}.
The branch complex power flow is formulated in \eqref{eq_power_acw}, to \eqref{eq_power_acw_q_2}, which depends on the voltage at each end of a branch. Note that problem induced by \eqref{AC_W_PF}, that is when $\mathcal{X}_t \leftarrow \mathcal{X}_{AC,t}$ in \eqref{OPF}, better represents reality, but is a non-linear and non-convex optimization problem. {It is relevant to mention that, following industry practice, the objective function only considers the active power in the assessment of fuel costs. On the other hand, the reactive power is indirectly accounted for through its relationship with the active power and all other variables and constraints to reach an AC feasible (implementable) operating point.}

The AC-OPF model is a non-convex non-linear problem (NLP), not suitable for the classical SDDP algorithm. Thus, as in many applications, convex approximations and relaxations can be used to meet the SDDP convexity hypothesis \cite{molzahn2019survey}. 
In general, those formulations are simplifications of the full AC model, and each one of them focuses on some particularities of the original problem. As a result, it is relevant to understand the tradeoff between each approximation quality and model tractability.

\textit{Relaxations} of the non-linear power flow constraints, when solved to optimality, provide valid bounds to the original problem because their feasible sets include all the solutions of the original problem. Convex relaxations are especially useful because their solutions are globally optimal for the relaxed problem, and the cuts generated by these relaxations are valid outer approximations for the real problem. However, as outer approximations, these cuts might lead to optimistic solutions. Although many convex relaxations exist for the optimal power flow problem, we focus on a limited subset. One simple linear relaxation used in the Brazilian official dispatch tools is the NFA, in which power flow limits are considered for each line, but KVL is ignored. 
In this lossless model, all the energy that is injected from an arbitrary bus $n$ into a line $(n,m)$, outputs at the receiving bus $m$. In other words, we make $\mathcal{X}_t \leftarrow \mathcal{X}_{NFA,t}$ in \eqref{OPF} where:
\begin{subequations} \label{NFA_PF}
    \begin{align}
        \mathcal{X}_{NFA,t} = \Big\{ \mathbf{x}_t \Big\vert \, f_{(n,m)t} = - f_{(m,n)t}\;\; \forall (n,m) \in \mathcal{L}\Big\}
    \end{align}
\end{subequations}

One of the more sophisticated convex relaxations is the \textit{semidefinite programming relaxation} (SDP). This formulation is obtained by replacing \eqref{eq_ac_w} by $\mathbf{W}_t \succeq 0 \land rank(\mathbf{W}_t)=1$ and then dropping the rank constraint, which is responsible for the non-convexity \cite{BAI2008383}. Hence, we define:
\begin{subequations} \label{SDP_PF}
    \begin{align}
        \mathcal{X}_{SDP,t}= \Big\{ \mathbf{x}_t  \Big\vert \, \exists \; q_t,\mathbf{W}_t,f^q_t : \text{\eqref{eq_voltage_acw_bounds}--\eqref{eq_gen_q_bounds_acw}}, \mathbf{W}_t \succeq 0  \Big\} .
    \end{align}
\end{subequations}

The \textit{Second-Order Cone relaxation} (SOC) is a non-linear convex relaxation that is tighter than NFA, but looser than SDP, \cite{1664986}. The feasible region of the AC formulation, shown in \eqref{AC_W_PF}, is contained within the feasible region of the SDP relaxation \eqref{SDP_PF}, which is contained in the SOC set \cite{low2014convex,low2014convex2}.
This time, \eqref{eq_ac_w} is replaced by: $|\mathrm{w}_{nmt}|^2 \leq  |v_{nt}|^2|v_{mt}|^2 = \mathrm{w}_{nnt}\mathrm{w}_{mmt}$.
The resulting problem may be specified as a second-order cone formulation, shown in \eqref{SOC_PF}.\\

$\mathcal{X}_{SOC,t} = \Big\{ \mathbf{x}_t \Big\vert \, \exists \; q_t,\mathbf{W}_t,f^q_t$: \eqref{eq_voltage_acw_bounds}--\eqref{eq_gen_q_bounds_acw},
\begin{subequations} \label{SOC_PF}
    \begin{align}
        & |\mathrm{w}_{nmt}|^2 \leq  \mathrm{w}_{nnt}\mathrm{w}_{mmt} \;\; \forall(n,m) \in \mathcal{L} \Big\} .
    \end{align}
\end{subequations}

Instead of relaxing, it is possible to approximate the non-linear Kirchhoff's Voltage Law equations through $\mathcal{X}_{DC,t}$, a linear DC power flow \cite{Stott2009_DC}. Alternatively, $\mathcal{X}_{DC,t}$ can be used to consider the quadratic DCLL approximation, which provides more accurate results \cite{coffrin2012approximating}, and can also be implemented through piecewise linear approximations. Both implementations are standard in power system and for the sake of brevity are omitted. We refer to \cite{Coffrin2018PowerModels} for the practical implementations of these models.

\vspace{-0.2cm}
\subsection{Open-source Julia package: HydroPowerModels.jl}

Solving a hydrothermal dispatch problem depends on the SDDP algorithm. However, until recently, there was no fast, reliable, and open-source implementation of the SDDP algorithm for the LTHDP. Without such a tool, researchers and practitioners have not had a common ground for discussing and analyzing different hydrothermal dispatch formulations and their solutions. Together with this work we made available an open-source tool, called \href{https://github.com/andrewrosemberg/HydroPowerModels.jl}{\textit{HydroPowerModels.jl}} \cite{Rosemberg2020}, that can be this common ground. \textit{HydroPowerModels.jl} can be used to assess the impact of modeling choices during the planning of a hydrothermal power system, as we show in the next section. To help time-inconsistency analysis, we had to extend the code developed in \textit{HydroPowerModels.jl} \cite{Rosemberg2020} and \textit{SDDP.jl} \cite{dowson_sddp.jl} significantly to implement the \textit{Fast Algorithm} proposed in \cite{7859406}.


\section{Case Studies}
\label{sec:cases}

In this section, the quality of the five approximations for the network constraints under study are compared. The time-inconsistency gap and other relevant operative indexes such as reservoir levels, thermal generation, and spot prices are studied to provide a more in-depth understanding of their differences. 

The inconsistent policies will be denoted according to the pair of network models used in the planning and implementation phases. In this work, all policies are compared with the AC power-flow model \eqref{AC_W_PF} using the same set of {3000} out-of-sample scenarios (scenarios not used in the planning stage where the cost-to-go functions were estimated).
Thus, for instance, the policy that used the transportation NFA model in the planning phase is named \textbf{NFA-AC inconsistent policy}. Analogously, the other policies will be named: \textbf{SOC-AC inconsistent policy}, \textbf{SDP-AC inconsistent policy}, \textbf{DC-AC inconsistent policy}, and \textbf{DCLL-AC inconsistent policy}. Finally, to estimate the time-inconsistency gap, we will also assess the cost for the planning policies, namely, \textbf{NFA planning policy}, \textbf{SOC planning policy}, \textbf{SDP planning policy}, the \textbf{DC planning policy}, and \textbf{DCLL planning policy}, each of which relies solely on their respective relaxations for both planning and implementation. 

{Although the proposed comparison method and metrics are general enough to embrace different applications, in this work we focus on high-voltage networks with higher X/R (reactance/resistance) rates in comparison to distribution networks. Therefore, the conclusions drawn in this work are conditioned to this setting. Additionally, because we are addressing an empirical study, the results and insights obtained in this work can only be ensured for the specific cases analysed in this work. Notwithstanding, we bring two examples of networks that showcase how different network features impact on the performance of each network model.}

\vspace{-0.1cm}
\subsection{Three-bus test system with loop}

As a case study, we use the three-bus system from \cite{7859406} to illustrate the effects of the underlying policies. A definition of case parameters following our notation can be found in \cite{AndrewThesis}. All three buses are connected in a loop, {thereby constituting an interesting example for evidencing the effect of KVL constraints on the quality of each studied approximation}. One hydro unit is located at bus $1$, one {thermoelectric} unit at bus $2$, and the most expensive thermoelectric unit and the demand at bus $3$. The planning horizon is 48 periods and the number of hours at each stage is 730 (one month). We use 3 scenarios per stage (low, medium and high), with similar values to \cite{7859406}, and simulate a single scenario {out of the $3^{48}$} per iteration in the forward step of the SDDP procedure.

Table \ref{table_costs_case3} shows the following information in its five columns: names of the inconsistent policies; the expected cost of the planning policies; the expected cost of the inconsistent policies at the implementation step (with the out-of-sample scenarios and AC power flow); the time-inconsistency GAP (the difference between the implementation and planning costs); and finally, the total computing times took for converging the SDDP with each network approximation (which is a measure of the computational burden that would be faced by system operators opting a given network representation).
\begin{table}[!ht]
\small
\renewcommand{\arraystretch}{1.3}
\centering
\caption{3-Bus Case: Planning and Implementation Policy Comparison.}
\label{table_costs_case3}
\begin{tabular}{c c c c c}
\hline
 \multicolumn{1}{p{0cm}}{\centering Policy \\ (plan,imp)} & \multicolumn{1}{p{0cm}}{\centering Planning \\ ($10^6 \$$)} & \multicolumn{1}{p{2cm}}{\centering Implementation \\ ($10^6 \$$)} & \multicolumn{1}{p{0cm}}{\centering GAP \\ (\%)} & \multicolumn{1}{p{0cm}}{\centering Time \\ (min.)} \\
\hline 
NFA-AC & $43.2317$ & $54.7234$ & $26.4$  & $5.40$\\
SOC-AC & $45.0477$ & $54.9197$ & $21.7$  & $23.92$\\
SDP-AC & $45.4524$ & $48.2853$  & $6.1$  & $47.47$\\
DC-AC & $43.3940$ & $48.8891$  & $12.5$  & $5.41$\\
DCLL-AC & $45.4158$ & $48.0700$ & $5.68$  & $16.09$\\
\hline
\end{tabular}
\end{table}

\begin{figure*}[htb]
\centerline{\includegraphics[width=0.85\linewidth]{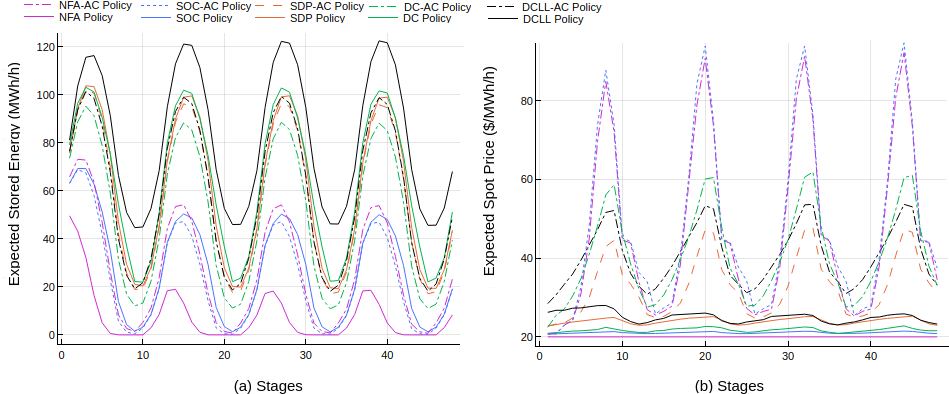}}
\caption{Three-bus case: (a) Expected Reservoir levels, (b) expected spot prices}
\label{fig:exp_case3}
\end{figure*}

The first important remark from Table \ref{table_costs_case3} is that the GAP increases with the simplification (relaxation) level. This indicates that the more optimistic one is in the planning phase, the higher will be the expected regret, i.e., the deviation of the implementation cost with regard to the planning one.
The inclusion relation between relaxations are the following: NFA is a relaxation of all the others, SOC is a relaxation of the SDP, and DC is a relaxation of the DCLL. 

The NFA-AC inconsistent policy, the simplest planning model, also features the lowest time, which is the main reason why the NFA simplification is widely adopted in practical studies \cite{maceiral2018twenty}. This simplification induces high operative costs in the implementation phase when an AC power flow is used to simulate the policy. 

The SOC-AC inconsistent policy features the second-highest GAP. Despite convex, this formulation is not easy to solve having the second higher computing time. This 3-Bus case is fully connected, hence it is one of the worst types of network for the SOC relaxation, opposed to radial systems \cite{1664986}. Interestingly, the SOC-AC policy performs slightly worse than the NFA-AC policy in the implementation phase, despite having a higher cost at the planning phase, but the difference is not statistically significant.
This suggests that the SOC approximation might exhibit inferior performances in the presence of very meshed grids {if no additional enhancements are considered. For instance, in \cite{kocuk2016strong}, a series of valid constraints are proposed and studied to improve the quality of SOC relaxations towards the tighter SDP model. Although it is beyond the scope of this paper to explore further variants of the selected approximations, it is worth highlighting that additional improvements in the SOC formulation would increase its complexity and thereby its computational burden, while not improving the policy quality (implementation cost and inconsistency gap) beyond that obtained with the SDP relaxations.}

The SDP-AC inconsistent policy has the tightest convex relaxation and the second-lowest GAP and cost in the implementation phase. Also, the implementation cost is almost equal to the lowest implementation cost (DCLL-AC). This benefit comes with the cost of the highest computational burden. The SDP solution tractability issues are due to appropriate solvers not being as evolved as linear and quadratic commercial solvers at the time of this study. 
The {DC-AC inconsistent policy} performs better than the SOC-AC and worse than the SDP-AC in terms of expected implementation cost. Moreover, it is a competitive policy having a GAP less than $2\%$ higher than the DCLL-AC policy. Although it does not model transmission losses,
it is capable of capturing the relevant operative constraints imposed by KVL and obtaining opportunity costs similar to those of the SDP model, but faster. 
The {DCLL-AC inconsistent policy} results indicate the best performance in terms of implementation cost, GAP, and average computing times for a small system. It puts together the DC model's capability to approximate KVL constraints and a reasonable description of the transmission losses. These results demonstrate a quadratic approximation can perform better than the conic SDP relaxation. Additionally, DCLL is more than two times faster than SDP, thereby representing an interesting alternative in the presence of meshed grids with non-negligible losses.

To further analyze these results, the expected storage is depicted in Figure \ref{fig:exp_case3} (a)
and the expected spot prices in Figure \ref{fig:exp_case3} (b) for all analyzed policies. 
In Figure \ref{fig:exp_case3} (a), notice that the {NFA planning policy} is the most optimistic (relaxed) model, leading to an aggressive use of the water that can not be implemented in practice due to the KVL constraints. This produces the highest inconsistency gap between the planning and the implemented policies depicted in both reservoir levels and spot prices. Conversely, the SOC model depicts an apparently consistent policy in terms of reservoir level. However, this relaxation precisely affects the KVL representation in the model, by further relaxing (2.a) dropping the semidefinite parcel of this constraint. Hence, in the presence of cycles, this model also exhibits a myopic view of KVL constraints, driving the system to extremely low expected reservoirs levels as the NFA does. Based on these two relaxations, the system operator, without acknowledging the electric constraints in the future, exposes the system to dangerous operating points and infeasible dispatches (requiring load curtailment). These infeasibilities cause a significant increase in the average spot prices as shown in Figure \ref{fig:exp_case3} (b). This relevant market distortion in spot prices, which is also accompanied by similar distortions in thermal generation, highlights the risk of ignoring KVL constraints in the planning stage.

In contrast to the aforementioned formulations ignoring the KVL constraints, there are the DC, SDP, and DCLL models. All of them consider approximations of the KVL constraints, thereby being aware of the operative difficulties of future stages caused by this constraints. Consequently, in the planning stage, their policies value the water accordingly, saving significantly more water than the KVL-myopic policies (NFA and SOC) as shown in \ref{fig:exp_case3} (a). This better representation of systems' constraints drives the reservoirs to safer levels enabling the system operator to circumvent the discrepancies between the planning and implementation models. As a consequence, spot prices distortions are significantly mitigated. The {DCLL planning policy} is even more pessimistic than the SDP and stores more water than the DC model. Still, the more accurate representation of the SDP relaxation allows a more efficient use of the stored water, leading to even lower discrepancies between planned and implemented spot prices.

According to Figure \ref{fig:exp_case3} (b), significant structural differences are found in the average spot prices between the respective planning and implementation policies. This stems from the fact that additional and expensive dispatches are needed to compensate for the optimistic view of the approximations in critical states (low reservoir levels). In the same Figure \ref{fig:exp_case3} (b), the {NFA planning policy} presents the lowest prices, as expected for the most relaxed problem that uses water resources as if no electrical constraints exist. The {SOC planning policy} has the second-lowest price given that it does not provide accurate representations of the network in the presence of cycles. In the sequel, the {DC planning policy}, {SDP planning policy}, and {DCLL planning policy} still provide a simplified version of the true network, albeit their representation are capable of considerably reducing the spot-price spikes {when implementing the policy under the more accurate AC power flow}. These better behaved spot-price profiles {observed in the SDP-AC and DCLL-AC implementation policies} stem from the better representation of the electrical constraints {in the planning phase. In other words, although increasing the cost in the planning phase, the better representation of the network allows the system to achieve better states and operating points in the implementation phase, resulting in lower costs and reducing market distortions}. Furthermore, the more inconsistent is the policy, the higher is the chance we will see spot-price spikes. This can be seen by comparing the {SDP-AC inconsistent policy}, the {SOC-AC inconsistent policy}, and the {NFA-AC inconsistent policy} which, in this order, incrementally relax the electric constraints in the planning stage. Notice that the {SOC planning policy} and the {SOC-AC inconsistent policy} differ here since the electrical operation provided by the planning policy is infeasible even though it has found an implementable storage management schedule on average. 
\vspace{-0.2cm}
\subsection{28-bus case study}
In order to further analyze the impacts of time-inconsistency due to network formulations {and the scalability of different network models within the SDDP technique}, we now use a larger case study using realistic data from the Bolivian system with the following characteristics: 28 buses, 26 loads, 34 generators (11 hydro generators), and 31 branches. The system is mostly radial with only 3 cycles. The planning horizon is 96 periods and we use 165 scenarios per stage derived from past data, and stagewise independently simulated in the forward iteration using a single scenario per iteration of the SDDP procedure in steps 1 and 2 of the algorithm.

\begin{table}[!ht]
\renewcommand{\arraystretch}{1.3}
\centering
\caption{28-bus case: Planning and Implementation Policy Comparison.}
\label{table_costs_bolivia}
\begin{tabular}{c c c c c}
\hline
 \multicolumn{1}{p{0cm}}{\centering Policy \\ (plan,imp)} & \multicolumn{1}{p{0cm}}{\centering Planning \\ ($10^6 \$$)} & \multicolumn{1}{p{2cm}}{\centering Implementation \\ ($10^6 \$$)} & \multicolumn{1}{p{0cm}}{\centering GAP \\ (\%)} & \multicolumn{1}{p{0cm}}{\centering Time \\ (min.)} \\
\hline  
NFA-AC & $0.58687$ & $0.68735$ & $17.067$ & $4.67$ \\
SOC-AC & $0.60818$ & $0.60935$ & $0.193$ & $119.23$ \\
DC-AC & $0.58714$ & $0.64523$ & $9.855$ & $5.22$ \\
DCLL-AC & $0.61310$ & $0.60941$ & $-0.602$ & $54.33$ \\
\hline
\end{tabular}
\end{table}

The {NFA-AC inconsistent policy} is still responsible for the higher expected cost in the implementation step when evaluated with the AC power flow. This table shows that the SOC convex relaxation produces a very small inconsistency GAP and one of the lowest operating costs. Under mild conditions, this relaxation is tight for radial systems (see \cite{1664986}). Thus, for the specific cases where cycles do not constrain much the least cost dispatch, which { differently from the previously studied system is the case here}, this relaxation performs reasonably well in terms of cost (policy quality). Nevertheless, this relaxation imposes the highest computational burden. It is worth mentioning that the SDP formulation has shown to be intractable for this problem and the hardware available. {The execution was interrupted after four weeks running without converging to reasonable gap values. Indeed, the SDP is known as one of ``the most difficult" classes of convex problems, and in spite of the relevant efforts, its scalability is still far behind what is possible with linear and quadratic programming. Furthermore, the computational burden of SDP-based relaxations for the AC power flow have been reported in \cite{kocuk2016strong}. Therefore, as the SDDP technique relies on the solution of a few million OPF problems,} the tractability of SDP-based policies poses a clear obstacle for the practical use of this network model, even for medium-sized systems.

\begin{figure*}[htb]
\centerline{\includegraphics[width=0.85\linewidth]{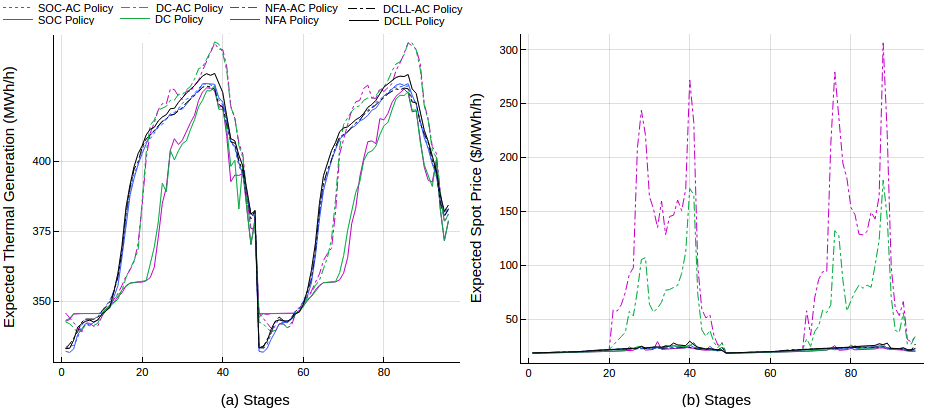}}
\caption{Case28: (a) expected thermal generation; (b) expected spot prices with time-inconsistency due to simplifications in
transmission-line modeling.}
\label{fig:exp_bolivia}
\end{figure*}

The {DC-AC inconsistent policy} performed better than the NFA-AC and worse than the {SOC-AC inconsistent policy} in terms of GAP and expected implementation cost. It also exhibits a significantly lower computational burden than the SOC relaxation and performance close to NFA. Thus, the DC network model provides an acceptable choice for practical applications with medium-sized systems. On the other hand, the {DCLL-AC inconsistent policy} features a statistically indistinguishable cost difference with the implementation policy based on the SOC relaxation (significance level of 0.05). Therefore, although the DCLL features a higher computational burden than the DC, it exhibits high quality solutions, the lowest inconsistency GAP, and a lower computational burden than the SOC relaxation (two times faster). It puts together the DC model capability of approximating KVL constraints and a reasonable description of the transmission losses. These results showcase that the DCLL approximation can have a close performance to the nonlinear SOC relaxation even in approximately radial systems.


To further analyze these results, the expected thermal generation is depicted in Figure \ref{fig:exp_bolivia} (a) and the expected spot prices in Figure \ref{fig:exp_bolivia} (b) for all analyzed policies. For NFA and DC, we see high thermal dispatches when evaluated under a AC power flow ({NFA-AC inconsistent policy} and {DC-AC inconsistent policy}) in Figure \ref{fig:exp_bolivia} (a) and, thus, also higher nodal prices are induced as seen in Figure \ref{fig:exp_bolivia} (b). { Differently from the previous case study, in this case, SOC and DCLL formulations provide similar performances in terms of inconsistency gap and market distortions. This is due to the system characteristics, namely, lower influence of KVL constrains on the dispatch decisions}. Therefore, their planned and implemented dispatches are very close to each other, mitigating the market distortions observed for the NFA and DC-based policies.

\vspace{-0.1cm}
\section{Conclusion}
\label{sec:conclusion}
In this work, we estimate and analyze the cost and impact of network simplifications in hydrothermal operation planning problems. We generalize and improve the results of \cite{7859406} by 1) considering the AC power flow to assess the inconsistency cost due to modeling simplifications, 2) considering five different network formulations in the cost assessment study, and 3) providing a novel open-source package \textit{HydroPowerModels.jl} \cite{Rosemberg2020}. For the first time in the literature, the quality of long-term hydrothermal dispatch policies based on conic network relaxations was evaluated and compared with policies based on DC models under the same basis, i.e., using the same out-of-sample scenarios and realistic AC power flow model. The results extend previously reported works showing that the optimistic assessments of widely adopted loss-less DC approximations also produce high dispatch costs and market distortions such as high and unexpected spot-price spikes and thermoelectric dispatches when evaluated under the more realistic AC power flow model. In this paper, we further identified the risk of poorly approximating the KVL constraints in the planning stage (when the opportunity cost of the water is calculated). We also find that the tightest network relaxation model (based on semidefinite programming) is not the one exhibiting the best operational performance. Instead, results show that the DC with quadratic line losses approximation exhibits the lowest system cost and inconsistency gaps.

{ Within the limitations of our case study, the results of our computational experiments allow us to draw the following conclusions:} 
\begin{itemize}
    \item The \textbf{Network Flow Approximation (NFA)} presents, on average, the lowest computational burden. However, this comes with the cost of exposing the system operator to high regrets due to the null representation of KVL constraints and losses. As a result, we find high operational costs in the actually implemented dispatches (frequently related to expensive and polluting thermal resources), and unjustified volatile spot-price profiles.
    
    \item The \textbf{Second Order Cone relaxation (SOC)} also provides poor approximations of the Kirchhoff's Voltage Law. Therefore, in the presence of cycles, this model also exhibits the typical distortions in spot prices and thermal generation dispatches. Its computational burden is higher than linear formulations, but significantly lower than the semidefinite relaxation. {In \cite{kocuk2016strong}, a series of valid constraints are proposed and studied to improve SOC relaxation results towards the results of the SDP relaxation. Although the consideration of such improvements would not change the conclusion that the DCLL provided the best tradeoff between computational burden and cost, we highlight the study of new methods and valid constraints to solve conic relaxations as a promising avenue for future research.}
    
    \item The \textbf{Semidefinite relaxation (SDP)} exhibits low distortions in both prices and thermal generation dispatches. However, this benefit comes with a high computational burden, which prevents its application in larger cases. Additionally, it is relevant to mention that although constituting the tightest relaxation, it is still a relaxation, thereby providing optimistic water values.
    
    \item The \textbf{DC approximation (DC)} is a standard and fair approach, performing reasonably well (small gap and low distortions in prices and dispatches) in all instances and with a reduced computational burden. 
    
    \item The \textbf{DC with line losses approximation (DCLL)} has the best performance in terms of implementation cost, it consistently presents the smallest inconsistency gaps, and exhibits very low distortions in prices and dispatches in comparison to the alternatives. The computational burden is not as reduced as in the DC case, as it is based on quadratic (convex) programming, but it is still faster than the two convex conic relaxations. Consequently, this approximation appears as a relevant alternative for the harder to solve conic relaxations and can bring significant benefits compared to the widely used Network Flow Approximation model currently adopted in Brazil. 
\end{itemize}
 
Based on these findings, we recommend system operators adopting network simplifications in the water value assessment to conduct further studies based on the DCLL network model under official models and data. As per our finds, we highlight the relevance of using the AC power-flow model and out-of-sample scenarios to conduct the long-term performance analyses.


\ifCLASSOPTIONcaptionsoff
  \newpage
\fi


\vspace{-0.3cm}

\bibliographystyle{IEEEtran}

\bibliography{ref.bib}


\end{document}